\documentclass{nagpreport}
\usepackage[british]{babel}

\usepackage{type1cm}

\usepackage{courier}

\usepackage[dvips]{graphicx}
\usepackage[dvips]{psfrag}
\usepackage[latin1]{inputenc}
\usepackage[T1]{fontenc}
\usepackage{mathptmx} 
\usepackage{appendix,url}       
\usepackage{float}  
\usepackage{amsmath,amsfonts,amssymb}   
\usepackage{easybmat}
\usepackage[small,bf,up]{caption} 
\usepackage[final]{microtype} 
\usepackage{array}
\usepackage{cite}      

\appendixtocon
\appendixtitleon
\appendixtitletocon
\appendixheaderon

\DeclareMathOperator*{\argmin}{arg\,min}
\DeclareMathOperator{\I}{I}

\begin{document}

\repno{07/18}

\newcommand{\tr}[1]{\ensuremath{\textnormal{trace}\left[{#1}\right]}}
\renewcommand{\trace}{\ensuremath{\textnormal{trace}}}
\newcommand{\veco}[1]{\ensuremath{\textnormal{vec}\left[{#1}\right]}}
\newcommand{\by}[2]{\ensuremath{\left[{#1}\times{#2}\right]}}
\newcommand{\kron}[2]{\ensuremath{\left({#1}\otimes{#2}\right)}}
\newcommand{\kronnp}[2]{\ensuremath{{#1}\otimes{#2}}} 
\renewcommand{\norm}[1]{\ensuremath{\left|\left|#1\right|\right|}}

\title{Kalman Filtering with Equality and Inequality State Constraints}

\author{Nachi Gupta\email{nachi} \and Raphael Hauser
\affil{Oxford University Computing Laboratory, Numerical Analysis Group,\\ Wolfson Building, Parks Road, Oxford OX1 3QD, U.K.}}

\kw{Constrained optimization, Kalman filtering, Nonlinear filters, Optimization methods, Quadratic programming, State estimation}

\ackn{The first author would like to thank the Clarendon Bursary for financial support.}


\summary{Both constrained and unconstrained optimization problems regularly appear in recursive tracking problems engineers currently address -- however, constraints are rarely exploited for these applications.  We define the Kalman Filter and discuss two different approaches to incorporating constraints.  Each of these approaches are first applied to equality constraints and then extended to inequality constraints.  We discuss methods for dealing with nonlinear constraints and for constraining the state prediction.  Finally, some experiments are provided to indicate the usefulness of such methods.}

\maketitle


\cleardoublepage

\section{Introduction}

Kalman Filtering \cite{Kalman1960} is a method to make real-time predictions for systems with some known dynamics.  Traditionally, problems requiring Kalman Filtering have been complex and nonlinear.  Many advances have been made in the direction of dealing with nonlinearities (e.g., Extended Kalman Filter \cite{BLK2001}, Unscented Kalman Filter \cite{JU1997}). These problems also tend to have inherent state space {\em equality\/} constraints (e.g., a fixed speed for a robotic arm) and state space {\em inequality\/} constraints (e.g., maximum attainable speed of a motor).  In the past, less interest has been generated towards constrained Kalman Filtering, partly because constraints can be difficult to model.  As a result, constraints are often neglected in standard Kalman Filtering applications.

The extension to Kalman Filtering with known equality constraints on the state space is discussed in \cite{SAP1988, TS1988, SC2002, WCC2002, Gupta2007}.  In this paper, we discuss two distinct methods to incorporate constraints into a Kalman Filter.  Initially, we discuss these in the framework of equality constraints.  The first method, projecting the updated state estimate onto the constrained region, appears with some discussion in \cite{SC2002,Gupta2007}.  We propose another method, which is to restrict the optimal Kalman Gain so the updated state estimate will not violate the constraint.  With some algebraic manipulation, the second method is shown to be a special case of the first method.

We extend both of these concepts to Kalman Filtering with inequality constraints in the state space.  This generalization for the first approach was discussed in \cite{SS2005}.\footnote{The similar extension for the method of \cite{WCC2002} was made in \cite{GHJ2005}.}   Constraining the optimal Kalman Gain was briefly discussed in \cite{Q1989}.  Further, we will also make the extension to incorporating state space constraints in Kalman Filter predictions.


Analogous to the way a Kalman Filter can be extended to solve problems containing non-linearities in the dynamics using an Extended Kalman Filter by linearizing locally (or by using an Unscented Kalman Filter), linear inequality constrained filtering can similarly be extended to problems with nonlinear constraints by linearizing locally (or by way of another scheme like an Unscented Kalman Filter).  The accuracy achieved by methods dealing with nonlinear constraints will naturally depend on the structure and curvature of the nonlinear function itself.  In the two experiments we provide, we look at incorporating inequality constraints to a tracking problem with nonlinear dynamics.

\section{Kalman Filter} \label{sec::kf}

A discrete-time Kalman Filter \cite{Kalman1960} attempts to find the best running estimate for a recursive system governed by the following model\footnote{The subscript $k$ on a variable stands for the $k$-th time step, the mathematical notation $\mathcal{N}\left(\mu,\Sigma\right)$ denotes a normally distributed random vector with mean $\mu$ and covariance $\Sigma$, and all vectors in this paper are column vectors (unless we are explicitly taking the transpose of the vector).}:

\begin{equation} \label{kfsm} 
x_{k} = F_{k,k-1} x_{k-1} + u_{k,k-1}, \qquad u_{k,k-1} \sim \mathcal{N}\left(0,Q_{k,k-1}\right) 
\end{equation}

\begin{equation} \label{kfmm} 
z_{k} = H_{k} x_{k} + v_{k}, \qquad v_{k} \sim \mathcal{N}\left(0,R_{k}\right) 
\end{equation}

Here $x_{k}$ is an $n$-vector that represents the true state of the underlying system and $F_{k,k-1}$ is an $n \times n$ matrix that describes the transition dynamics of the system from  $x_{k-1}$ to $x_{k}$.  The measurement made by the observer is an $m$-vector $z_{k}$, and $H_{k}$ is an $m \times n$ matrix that transforms a vector from the state space into the appropriate vector in the measurement space.  The noise terms $u_{k,k-1}$ (an $n$-vector) and $v_{k}$ (an $m$-vector) encompass known and unknown errors in $F_{k,k-1}$ and $H_{k}$ and are normally distributed with mean 0 and covariances given by $n \times n$ matrix $Q_{k,k-1}$ and $m \times m$ matrix $R_{k}$, respectively.  At each iteration, the Kalman Filter makes a state prediction for $x_k$, denoted $\hat{x}_{k|k-1}$.  We use the notation ${k|k-1}$ since we will only use measurements provided until time-step $k-1$ in order to make the prediction at time-step $k$.  The state prediction error $\tilde{x}_{k|k-1}$ is defined as the difference between the true state and the state prediction, as below.

\begin{equation} \label{se1}
\tilde{x}_{k|k-1} = x_{k} - \hat{x}_{k|k-1}
\end{equation}

The covariance structure for the expected error on the state prediction is defined as the expectation of the outer product of the state prediction error.  We call this covariance structure the error covariance prediction and denote it $P_{k|k-1}$.\footnote{We use the prime notation on a vector or a matrix to denote its transpose throughout this paper.}

\begin{equation} \label{P-outer1}
P_{k|k-1} = \mathbb{E}\left[\left(\tilde{x}_{k|k-1}\right)\left(\tilde{x}_{k|k-1}\right)'\right]
\end{equation}

The filter will also provide an updated state estimate for $x_{k}$, given all the measurements provided up to and including time step $k$.  We denote these estimates by $\hat{x}_{k|k}$.  We similarly define the state estimate error $\tilde{x}_{k|k}$ as below.

\begin{equation} \label{se2}
\tilde{x}_{k|k} = x_{k} - \hat{x}_{k|k}
\end{equation}

The expectation of the outer product of the state estimate error represents the covariance structure of the expected errors on the state estimate, which we call the updated error covariance and denote $P_{k|k}$.

\begin{equation} \label{P-outer2}
P_{k|k} = \mathbb{E}\left[\left(\tilde{x}_{k|k}\right)\left(\tilde{x}_{k|k}\right)'\right]
\end{equation}

At time-step $k$, we can make a prediction for the underlying state of the system by allowing the state to transition forward using our model for the dynamics and noting that $\mathbb{E}\left[u_{k,k-1}\right] = 0$.  This serves as our state prediction.

\begin{equation} \label{kfsp} 
\hat{x}_{k|k-1} = F_{k,k-1} \hat{x}_{k-1|k-1} 
\end{equation}

If we expand the expectation in Equation \eqref{P-outer1}, we have the following equation for the error covariance prediction.

\begin{equation} \label{kfcp} 
P_{k|k-1} = F_{k,k-1} P_{k-1|k-1} F_{k,k-1}' + Q_{k,k-1}
\end{equation}

We can transform our state prediction into the measurement space, which is a prediction for the measurement we now expect to observe.

\begin{equation} \label{kfmp} 
\hat{z}_{k|k-1} = H_{k} \hat{x}_{k|k-1}
\end{equation}

The difference between the observed measurement and our predicted measurement is the measurement residual, which we are hoping to minimize in this algorithm.

\begin{equation} \label{kfi} 
\nu_{k} = z_{k} - \hat{z}_{k|k-1} 
\end{equation}

We can also calculate the associated covariance for the measurement residual, which is the expectation of the outer product of the measurement residual with itself, $\mathbb{E}\left[\nu_k \nu_k'\right]$.  We call this the measurement residual covariance.

\begin{equation} \label{kfic} 
S_{k} = H_{k} P_{k|k-1} H_{k}' + R_{k} 
\end{equation}

We can now define our updated state estimate as our prediction plus some perturbation, which is given by a weighting factor times the measurement residual.  The weighting factor, called the Kalman Gain, will be discussed below.

\begin{equation} \label{kfsu} 
\hat{x}_{k|k} = \hat{x}_{k|k-1} + K_{k}  \nu_{k}
\end{equation}

Naturally, we can also calculate the updated error covariance by expanding the outer product in Equation \eqref{P-outer2}.\footnote{The $\I$ in Equation \eqref{kfcu} represents the $n \times n$ identity matrix.  Throughout this paper, we use $\I$ to denote the same matrix, except in Appendix~\ref{app::kv}, where $\I$ is the appropriately sized identity matrix.}

\begin{equation} \label{kfcu} 
P_{k|k} = \left(\I - K_{k} H_{k}\right) P_{k|k-1}   \left(\I - K_{k} H_{k}\right)' + K_k R_k K_k'
\end{equation}

Now we would like to find the Kalman Gain $K_k$, which minimizes the mean square state estimate error, $\mathbb{E}\left[\left|\tilde{x}_{k|k}\right|^2\right]$.  This is the same as minimizing the trace of the updated error covariance matrix above.\footnote{Note that $v'v = \tr{vv'}$ for some vector $v$.}  After some calculus, we find the optimal gain that achieves this, written below.\footnote{We could also minimize the mean square state estimate error in the $N$ norm, where $N$ is a positive definite and symmetric weighting matrix.  In the $N$ norm, the optimal gain would be $K^N_k = N^{\frac{1}{2}}K_k$.}

\begin{equation} \label{kfkg} 
K_{k} = P_{k|k-1} H_{k}' S_{k}^{-1} 
\end{equation}

The covariance matrices in the Kalman Filter provide us with a measure for uncertainty in our predictions and updated state estimate.  This is a very important feature for the various applications of filtering since we then know how much to trust our predictions and estimates.  Also, since the method is recursive, we need to provide an initial covariance that is large enough to contain the initial state to ensure comprehensible performance.  For a more detailed discussion of Kalman Filtering, we refer the reader to the following book \cite{BLK2001}.

\section{Equality Constrained Kalman Filtering}

A number of approaches have been proposed for solving the equality constrained Kalman Filtering problem \cite{TS1988, SAP1988, WCC2002, SC2002, Gupta2007}.  In this paper, we show two different methods.  The first method will restrict the state at each iteration to lie in the equality constrained space.  The second method will start with a constrained prediction, and restrict the Kalman Gain so that the estimate will lie in the constrained space.   Our equality constraints in this paper will be defined as below, where $A$ is a $q \times n$ matrix, $b$ a $q$-vector, and $x_k$, the state, is a $n$-vector.\footnote{$A$ and $b$ can be different for different $k$.  We don't subscript each $A$ and $b$ to avoid confusion.}

\begin{equation} \label{constraints}
A x_k = b
\end{equation}

So we would like our updated state estimate to satisfy the constraint at each iteration, as below.

\begin{equation} \label{kfsu-con}
A \hat{x}_{k|k} = b
\end{equation}

Similarly, we may also like the state prediction to be constrained, which would allow a better forecast for the system.

\begin{equation}
A \hat{x}_{k|k-1} = b
\end{equation}

In the following subsections, we will discuss methods for constraining the updated state estimate.  In Section~\ref{sec::aic}, we will extend these concepts and formulations to the inequality constrained case, and in Section~\ref{sec::csp}, we will address the problem of constraining the prediction, as well.

\subsection{Projecting the state to lie in the constrained space} \label{sec::pue}



We can solve the following minimization problem for a given time-step $k$, where $\hat{x}_{k|k}^{P}$ is the constrained estimate,  $W_k$ is any positive definite symmetric weighting matrix, and $\hat{x}_{k|k}$ is the unconstrained Kalman Filter updated estimate.
 
\begin{equation} \label{eq-proj-problem}
\hat{x}_{k|k}^{P} = \argmin_{x \in \mathbb{R}^n} \ \left\{\left(x - \hat{x}_{k|k} \right)' W_k \left(x - \hat{x}_{k|k} \right) : A x = b\right\}
\end{equation}

The best constrained estimate is then given by

\begin{equation} \label{bce-xP}
\hat{x}_{k|k}^{P} = \hat{x}_{k|k} - W_k^{-1} A' \left( A W_k^{-1} A' \right)^{-1} \left(A \hat{x}_{k|k} - b \right)
\end{equation}

To find the updated error covariance matrix of the equality constrained filter, we first define the matrix $\Upsilon$ below.\footnote{Note that $\Upsilon A$ is a projection matrix, as is $\left(\I - \Upsilon A\right)$, by definition.  If $A$ is poorly conditioned, we can use a QR factorization to avoid squaring the condition number.}

\begin{equation}
\Upsilon = W_k^{-1} A' \left(A W_k^{-1} A' \right)^{-1}
\end{equation}

Equation \eqref{bce-xP} can then be re-written as following.

\begin{equation} \label{xeq}
\hat{x}_{k|k}^P = \hat{x}_{k|k} - \Upsilon\left(A \hat{x}_{k|k} - b \right)
\end{equation}

We can find a reduced form for $x_k - \hat{x}_{k|k}^P$ as below.

\begin{subequations}
\begin{align}
x_k - \hat{x}_{k|k}^P &= x_k - \hat{x}_{k|k} +\Upsilon \left(A \hat{x}_{k|k} - b - \left(A x_k - b \right)\right) \\
&=
	x_k - \hat{x}_{k|k} +\Upsilon \left(A \hat{x}_{k|k} - A x_k\right) \\
&=
	-\left(\I - \Upsilon A \right) \left(\hat{x}_{k|k} - x_k\right)
\end{align}
\end{subequations}

Using the definition of the error covariance matrix, we arrive at the following expression.

\begin{subequations} \label{bce-PP}
\begin{align}
P_{k|k}^P &= \mathbb{E}\left[\left(x_k - \hat{x}_{k|k}^P\right)\left(x_k - \hat{x}_{k|k}^P\right)'\right] \\
&= 
	\mathbb{E}\left[\left(\I - \Upsilon A \right) \left(\hat{x}_{k|k} - x_k\right) \left(\hat{x}_{k|k} - x_k\right)' \left(\I - \Upsilon A \right)'\right] \\
&= 
	\left(\I - \Upsilon A \right) P_{k|k} \left(\I - \Upsilon A \right)' \\
&= 
	P_{k|k} - \Upsilon A P_{k|k} - P_{k|k} A' \Upsilon' +  \Upsilon A P_{k|k} A' \Upsilon' \\
&= 
	P_{k|k} - \Upsilon A P_{k|k} \\
&=  \label{Peq}
	\left(\I - \Upsilon A \right) P_{k|k} 
\end{align}
\end{subequations}

It can be shown that choosing $W_k = P_{k|k}^{-1}$ results in the smallest updated error covariance.  This also provides a measure of the information in the state at $k$.\footnote{If $M$ and $N$ are covariance matrices, we say $N$ is smaller than $M$ if $M-N$ is positive semidefinite.  Another formulation for incorporating equality constraints into a Kalman Filter is by observing the constraints as pseudo-measurements \cite{TS1988,WCC2002}.  When $W_k$ is chosen to be $P_{k|k}^{-1}$, both of these methods are mathematically equivalent \cite{Gupta2007}.  Also, a more numerically stable form of Equation \eqref{bce-PP} with discussion is provided in \cite{Gupta2007}.}

\subsection{Restricting the optimal Kalman Gain so the updated state estimate lies in the constrained space}

Alternatively, we can expand the updated state estimate term in Equation \eqref{kfsu-con} using Equation \eqref{kfsu}.  

\begin{equation} 
A \left(  \hat{x}_{k|k-1} + K_{k}  \nu_{k} \right) = b
\end{equation}

Then, we can choose a Kalman Gain $K_k^R$, that forces the updated state estimate to be in the constrained space.  In the unconstrained case, we chose the optimal Kalman Gain $K_k$, by solving the minimization problem below which yields Equation \eqref{kfkg}.

\begin{equation}
K_k = \argmin_{K \in \mathbb{R}^{n \times m}} \trace \left[ \left(\I - K H_{k}\right) P_{k|k-1}   \left(\I - K H_{k}\right)' + K R_k K'\right]
\end{equation}

Now we seek the optimal $K_k^R$ that satisfies the constrained optimization problem written below for a given time-step $k$.

\begin{equation} \label{min-con}
\begin{split}
 K_k^R = \argmin_{K \in \mathbb{R}^{n \times m}} & \trace \left[ \left(\I - K H_{k}\right) P_{k|k-1}   \left(\I - K H_{k}\right)' + K R_k K'\right] \\
 \textnormal{s.t. } & A \left(  \hat{x}_{k|k-1} + K  \nu_{k} \right) = b 
\end{split}
\end{equation}

We will solve this problem using the method of Lagrange Multipliers.  First, we take the steps below, using the vec notation (column stacking matrices so they appear as long vectors, see Appendix~\ref{app::kv}) to convert all appearances of $K$ in Equation \eqref{min-con} into long vectors.  Let us begin by expanding the following term.\footnote{Throughout this paper, a number in parentheses above an equals sign means we made use of this equation number.}




\begin{subequations}
\begin{gather}
\nonumber\trace\left[\left(\I - K H_{k}\right) P_{k|k-1}   \left(\I - K H_{k}\right)' + K R_k K' \right] \qquad \qquad \qquad \qquad \qquad \qquad \qquad\\
\begin{aligned}
&\stackrel{\hphantom{\eqref{kfic}}}{=}\trace \left[ P_{k|k-1}  - K H_{k} P_{k|k-1} - P_{k|k-1} H_{k}' K'  +  K H_{k} P_{k|k-1}  H_{k}' K' + K R_k K' \right] \\
&\stackrel{\eqref{kfic}}{=} \trace \left[ P_{k|k-1}  - K H_{k} P_{k|k-1} - P_{k|k-1} H_{k}' K' +  K S_k K' \right] \\
&\stackrel{\hphantom{\eqref{kfic}}}{=}\label{trace-separated}\trace \left[ P_{k|k-1} \right] - \trace \left[ K H_{k} P_{k|k-1} \right] - \trace \left[ P_{k|k-1} H_{k}' K' \right] + \trace \left[ K S_k K' \right]
\end{aligned}
\end{gather}
\end{subequations}

We now expand the last three terms in Equation \eqref{trace-separated} one at a time.\footnote{We use the symmetry of $P_{k|k-1}$ in Equation \eqref{KHP} and the symmetry of $S_k$ in Equation \eqref{KSK}.}

\begin{equation} \label{KHP}
\begin{aligned}
\trace \left[ K H_{k} P_{k|k-1} \right] 
\stackrel{\eqref{tr-ab}}{=} 
	\veco{\left(H_k P_{k|k-1}\right)'}' \veco{K} \\
=
	\veco{P_{k|k-1} H_k'}' \veco{K}
\end{aligned}
\end{equation}


\begin{equation}
\trace \left[ P_{k|k-1} H_{k}' K' \right]
\stackrel{\eqref{tr-ab}}{=} 
	\veco{K}' \veco{P_{k|k-1} H_k'}
\end{equation}





\begin{equation} \label{KSK}
\begin{aligned}
\trace \left[ K S_k K' \right] 
&\stackrel{\eqref{tr-ab}}{=} 
	\veco{K}' \veco{K S_k} \\
&\stackrel{\eqref{vec-ab}}{=} 
	\veco{K}' \kron{S}{\I} \veco{K}
\end{aligned}
\end{equation}

Remembering that $\trace \left[ P_{k|k-1} \right]$ is constant, our objective function can be written as below.

\begin{equation}
\begin{aligned}
\veco{K}' \left(\I \otimes S_k \right) \veco{K'} &- \veco{P_{k|k-1} H_k'}' \veco{K}\\
&- \veco{K}' \veco{P_{k|k-1} H_k'}
\end{aligned}
\end{equation}

Using Equation \eqref{vec-abc} on the equality constraints, our minimization problem is the following.



\begin{equation}
\begin{split}
K_k^R = \argmin_{K \in \mathbb{R}^{n \times m}}& \ \veco{K}' \kron{S_k}{\I} \veco{K} \\
&- \veco{P_{k|k-1} H_k'}' \veco{K} \\ 
& - \veco{K}'  \veco{P_{k|k-1} H_k'} \\
\textnormal{s.t. } &  \left( \nu_{k}' \otimes A \right) \veco{K}  = b - A \hat{x}_{k|k-1}
\end{split}
\end{equation}



Further, we simplify this problem so the minimization problem has only one quadratic term.  We complete the square as follows.  We want to find the unknown variable $\mu$ which will cancel the linear term.  Let the quadratic term appear as follows.  Note that the non-``$\veco{K}$" term is dropped as is is irrelevant for the minimization problem.


\begin{equation}
\left(\veco{K} + \mu \right)' \kron{S_k}{\I} \left( \veco{K} + \mu \right)
\end{equation}

The linear term in the expansion above is the following.

\begin{equation}
\veco{K}'  \kron{S_k}{\I} \mu + \mu' \kron{S_k}{\I} \veco{K}
\end{equation}

So we require that the two equations below hold.

\begin{equation}
\begin{aligned}
\kron{S_k}{\I} \mu &= -\veco{P_{k|k-1} H_k'} \\
\mu' \kron{S_k}{\I} &= -\veco{P_{k|k-1} H_k'}'
\end{aligned}
\end{equation}

This leads to the following value for $\mu$.

\begin{equation}
\begin{aligned}
\mu 
&\stackrel{\eqref{kron-inv}}{=}
	 - \kron{S_k^{-1}}{\I} \veco{P_{k|k-1} H_k'} \\
&\stackrel{\eqref{vec-abc}}{=}
	-\veco{P_{k|k-1} H_k' S_k^{-1}} \\
&\stackrel{\eqref{kfkg}}{=}
	-\veco{K_k}
\end{aligned}
\end{equation}





Using Equation \eqref{vec-sum}, our quadratic term in the minimization problem becomes the following.

\begin{equation}
\left(\veco{K - K_k} \right)' \kron{S_k}{\I} \left( \veco{K - K_k} \right)
\end{equation}

Let $l = \veco{K - K_k}$.  Then our minimization problem becomes the following.

\begin{equation}
\begin{aligned}
K_k^R = \argmin_{l \in \mathbb{R}^{mn}} & \ l' \kron{S_k}{\I} l \\
\textnormal{s.t. }&  \left( \nu_{k}' \otimes A \right) \left(l + \veco{K_{k}}\right)  = b - A \hat{x}_{k|k-1}
\end{aligned}
\end{equation}

We can then re-write the constraint taking the $\veco{K_k}$ term to the other side as below.

\begin{equation}
\begin{aligned}
\left( \nu_{k}' \otimes A \right) l & = b - A \hat{x}_{k|k-1} - \left( \nu_{k}' \otimes A \right) \veco{K_{k}} \\
& \stackrel{\eqref{vec-abc}}{=} b - A \hat{x}_{k|k-1} -\veco{A K_{k} \nu_k} \\
& = b - A \hat{x}_{k|k-1}  - A K_{k} \nu_k \\
& \stackrel{\eqref{kfsu}}= b - A \hat{x}_{k|k}
\end{aligned}
\end{equation}

This results in the following simplified form.

\begin{equation} \label{first-SDPT3}
\begin{aligned}
K_k^R = \argmin_{l \in \mathbb{R}^{mn}}&\ l' \kron{S_k}{\I} l \\
\textnormal{s.t. }&  \left( \nu_{k}' \otimes A \right) l  = b - A \hat{x}_{k|k}
\end{aligned}
\end{equation}

We form the Lagrangian $\mathcal{L}$, where we introduce $q$ Lagrange Multipliers in vector $ \lambda = \left( \lambda_1, \lambda_2, \ldots, \lambda_q \right)'$



\begin{equation}
\begin{aligned}
\mathcal{L} = & l' \kron{S_k}{\I} l -  \lambda' \left[\left( \nu_{k}' \otimes A \right) l  - b + A \hat{x}_{k|k}\right]
\end{aligned}
\end{equation}

We take the partial derivative with respect to $l$.\footnote{We used the symmetry of $\kron{S_k}{\I}$ here.}

\begin{equation} \label{partial1}
\frac{\partial \mathcal{L}}{\partial l} = 2 l' \kron{S_k}{\I} - \lambda' \left( \nu_{k}' \otimes A \right)  \\
\end{equation}

Similarly we can take the partial derivative with respect to the vector $\lambda$.

\begin{equation}
\frac{\partial \mathcal{L}}{\partial \lambda}  = \left( \nu_{k}' \otimes A \right) l  - b + A \hat{x}_{k|k}
\end{equation}

When both of these derivatives are set equal to the appropriate size zero vector, we have the solution to the system.  Taking the transpose of Equation \eqref{partial1}, we can write this system as $Mn = p$ with the following block definitions for $M,n$, and $p$.

\begin{equation} \label{M-matrix}
M = \begin{bmatrix}
	2  \kronnp{S_k}{\I} & \nu_{k} \otimes A' \\
	 \nu_{k}' \otimes A & 0_{\by{q}{q}}
\end{bmatrix}
\end{equation}

\begin{equation} \label{n-vector}
n = \begin{bmatrix}
	l \\
	\lambda
\end{bmatrix}
\end{equation}

\begin{equation} \label{p-vector}
p = \begin{bmatrix}
	0_{\by{mn}{1}} \\
	b - A \hat{x}_{k|k}
\end{bmatrix}
\end{equation}

We solve this system for vector $n$ in Appendix~\ref{app::Mnp}.  The solution for $l$ is pasted below.

\begin{equation}
\left(\left[S_k^{-1} \nu_k \left(\nu_{k}' S_k^{-1} \nu_k \right)^{-1}\right] \otimes \left[A' \left(A A' \right)^{-1} \right]\right) \left(b - A \hat{x}_{k|k}\right)
\end{equation}

Bearing in mind that $b - A \hat{x}_{k|k} = \veco{b - A \hat{x}_{k|k}}$, we can use Equation \eqref{vec-abc} to re-write $l$ as below.\footnote{Here we used the symmetry of $S_k^{-1}$ and $\left(\nu_{k}' S_k^{-1} \nu_k \right)^{-1}$ (the latter of which is actually just a scalar).}

\begin{equation}
\veco{A' \left(A A' \right)^{-1}\left(b - A \hat{x}_{k|k} \right)  \left(\nu_{k}' S_k^{-1} \nu_k \right)^{-1} \nu_k' S_k^{-1}}
\end{equation}

The resulting matrix inside the vec operation is then an $n$ by $m$ matrix.  Remembering the definition for $l$, we notice that $K - K_k$ results in an $n$ by $m$ matrix also.  Since both of the components inside the vec operation result in matrices of the same size, we can safely remove the vec operation from both sides.  This results in the following optimal constrained Kalman Gain $K_k^R$.

\begin{equation}
K_k - A' \left(A A' \right)^{-1}\left(A \hat{x}_{k|k} - b \right)  \left(\nu_{k}' S_k^{-1} \nu_k \right)^{-1} \nu_k' S_k^{-1}
\end{equation}

If we now substitute this Kalman Gain into Equation \eqref{kfsu} to find the constrained updated state estimate, we end up with the following.

\begin{equation}
\hat{x}_{k|k}^R = \hat{x}_{k|k} - A' \left(A A' \right)^{-1}\left(A \hat{x}_{k|k}  - b \right) 
\end{equation}

This is of course equivalent to the result of Equation \eqref{bce-xP} with the weighting matrix $W_k$ chosen as the identity matrix.  The error covariance for this estimate is given by Equation \eqref{bce-PP}.\footnote{We can use the unconstrained or constrained Kalman Gain to find this error covariance matrix.
Since the constrained Kalman Gain is suboptimal for the unconstrained problem, before projecting onto the constrained space, the constrained covariance will be different from the unconstrained covariance.  However, the difference lies exactly in the space orthogonal to which the covariance is projected onto by  Equation \eqref{bce-PP}.  The proof is omitted for brevity.}

\section{Adding Inequality Constraints} \label{sec::aic}
In the more general case of this problem, we may encounter equality and inequality constraints, as given below.\footnote{$C$ and $d$ can be different for different $k$.  We don't subscript each $C$ and $d$ to simplify notation.}


\begin{equation} \label{ineq-constraints}
\begin{split}
A x_{k} = b\\
C x_{k} \leq d 
\end{split}
\end{equation}

So we would like our updated state estimate to satisfy the constraint at each iteration, as below.

\begin{equation} 
\begin{split}
A \hat{x}_{k|k} = b \\
C \hat{x}_{k|k} \leq d
\end{split}
\end{equation}

Similarly, we may also like the state prediction to be constrained, which would allow a better forecast for the system.

\begin{equation}
\begin{split}
A \hat{x}_{k|k-1} = b \\
C \hat{x}_{k|k-1} \leq d
\end{split}
\end{equation}

We will present two analogous methods to those presented for the equality constrained case.  In the first method, we will run the unconstrained filter, and at each iteration constrain the updated state estimate to lie in the constrained space.  In the second method, we will find a Kalman Gain $\check{K}_k^R$ such that the the updated state estimate will be forced to lie in the constrained space.  In both methods, we will no longer be able to find an analytic solution as before.  Instead, we use numerical methods.

\subsection{By Projecting the Unconstrained Estimate} \label{sec::pue-ineq}

Given the best unconstrained estimate, we could solve the following minimization problem for a given time-step $k$, where $\check{x}_{k|k}^{P}$ is the inequality constrained estimate and $W_k$ is any positive definite symmetric weighting matrix.

\begin{equation}
\begin{aligned}
\check{x}_{k|k}^{P} =  \argmin_{x} &\  \left(x - \hat{x}_{k|k} \right)' W_k \left(x - \hat{x}_{k|k} \right) \\
\textnormal{s.t. } & A x = b \\
& C x \leq d
\end{aligned}
\end{equation}

For solving this inequality constrained optimization problem, we can use a variety of standard methods, or even an out-of-the-box solver, like \verb|fmincon| in Matlab.  Here we use an active set method \cite{Fletcher1981}.  This is a common method for dealing with inequality constraints, where we treat a subset of the constraints (called the active set) as additional equality constraints.  We ignore any inactive constraints when solving our optimization problem.  After solving the problem, we check if our solution lies in the space given by the inequality constraints.  If it doesn't, we start from the solution in our previous iteration and move in the direction of the new solution until we hit a set of constraints.  For each iteration, the active set is made up of those inequality constraints with non-zero Lagrange Multipliers.


We first find the best estimate (using Equation \eqref{bce-xP} for the equality constrained problem with the equality constraints given in Equation \eqref{ineq-constraints} plus the active set of inequality constraints.   Let us call the solution to this $\check{x}_{k|k,j}^{P*}$ since we have not yet checked if the solution lies in the inequality constrained space.\footnote{For the inequality constrained filter, we allow multiple iterations within each step.  The $j$ subscript indexes  these further iterations.}  In order to check this, we find the vector that we moved along to reach $\check{x}_{k|k,j}^{P*}$.  This is given by the following.

\begin{equation}
s = \check{x}_{k|k,j}^{P*} - \check{x}_{k|k,j-1}^P
\end{equation}

We now iterate through each of our inequality constraints, to check if they are satisfied.  If they are all satisfied, we choose $\tau_{\max}=1$. If they are not, we choose the largest value of $\tau_{\max}$ such that $\hat{x}_{k|k,j-1} + \tau_{\max} s$ lies in the inequality constrained space.  We choose our estimate to be

\begin{equation}
\check{x}_{k|k,j}^P = \check{x}_{k|k,j-1}^{P} + \tau_{\max} s
\end{equation}

If we find the solution has converged within a pre-specified error, or we have reached a pre-specified maximum number of iterations, we choose this as the updated state estimate to our inequality constrained problem, denoted $\check{x}_{k|k}^P$.  If we would like to take a further iteration on $j$, we check the Lagrange Multipliers at this new solution to determine the new active set.\footnote{The previous active set is not relevant.}  We then repeat by finding the best estimate for the equality constrained problem including the new active set as additional equality constraints.  Since this is a Quadratic Programming problem, each step of $j$ guarantees the same estimate or a better estimate.

When calculating the error covariance matrix for this estimate, we can also add on the safety term below.

\begin{equation}
\left(\check{x}_{k|k,j}^P - \check{x}_{k|k,j-1}^{P}\right)\left(\check{x}_{k|k,j}^P - \check{x}_{k|k,j-1}^{P}\right)'
\end{equation}

This is a measure of our convergence error and should typically be small relative to the unconstrained error covariance.  We can then use Equation \eqref{bce-PP} to project the covariance matrix onto the constrained subspace, but we only use the defined equality constraints.  We do not incorporate any constraints in the active set when computing Equation \eqref{bce-PP} since these still represent inequality constraints on the state.  Ideally we would project the error covariance matrix into the inequality constrained subspace, but this projection is not trivial.  


\subsection{By Restricting the Optimal Kalman Gain}











We could solve this problem by restricting the optimal Kalman gain also, as we did for equality constraints previously.  We seek the optimal $K_k$ that satisfies the constrained optimization problem written below for a given time-step $k$.  

\begin{equation} \label{min-con}
\begin{aligned}
\check{K}^R_k = \argmin_{K \in \mathbb{R}^{n \times m}} & \trace \left[\left(\I - K H_{k}\right) P_{k|k-1}   \left(\I - K H_{k}\right)' + K R_k K'\right] \\ 
\textnormal{s.t. } & A \left(  \hat{x}_{k|k-1} + K_{k}  \nu_{k} \right) = b \\
& C \left( \hat{x}_{k|k-1} + K_{k}  \nu_{k} \right) \leq d  
\end{aligned}
\end{equation}

Again, we can solve this problem using any inequality constrained optimization method (e.g., \verb|fmincon| in Matlab or the active set method used previously).  Here we solved the optimization problem using SDPT3, a Matlab package for solving semidefinite programming problems \cite{TTT1999}.  When calculating the covariance matrix for the inequality constrained estimate, we use the restricted Kalman Gain.  Again, we can add on the safety term for the convergence error, by taking the outer product of the difference between the updated state estimates calculated by the restricted Kalman Gain for the last two iterations of SDPT3.  This covariance matrix is then projected onto the subspace as in Equation \eqref{bce-PP} using the equality constraints only.




\section{Dealing with Nonlinearities} \label{sec::nl}

Thus far, in the Kalman Filter we have dealt with linear models and constraints.  A number of methods have been proposed to handle nonlinear models (e.g., Extended Kalman Filter \cite{BLK2001}, Unscented Kalman Filter \cite{JU1997}).  In this paper, we will focus on the most widely used of these, the Extended Kalman Filter.  Let's re-write the discrete unconstrained Kalman Filtering problem from Equations \eqref{kfsm} and \eqref{kfmm} below, incorporating nonlinear models.

\begin{equation} \label{kfsm-nl} 
x_{k} = f_{k,k-1} \left(x_{k-1}\right) + u_{k,k-1}, \qquad u_{k,k-1} \sim \mathcal{N}\left(0,Q_{k,k-1}\right) 
\end{equation}

\begin{equation} \label{kfmm-nl} 
z_{k} = h_{k} \left(x_{k}\right) + v_{k}, \qquad v_{k} \sim \mathcal{N}\left(0,R_{k}\right) 
\end{equation}

In the above equations, we see that the transition matrix $F_{k,k-1}$ has been replaced by the nonlinear vector-valued function $f_{k,k-1}\left(\cdot\right)$, and similarly, the matrix $H_k$, which transforms a vector from the state space into the measurement space, has been replaced by the nonlinear vector-valued function $h_k\left(\cdot\right)$.  The method proposed by the Extended Kalman Filter is to linearize the nonlinearities about the current state prediction (or estimate).  That is, we choose $F_{k,k-1}$ as the Jacobian of $f_{k,k-1}$ evaluated at $\hat{x}_{k-1|k-1}$, and $H_k$ as the Jacobian of $h_k$ evaluated at $\hat{x}_{k|k-1}$ and proceed as in the linear Kalman Filter of Section~\ref{sec::kf}.\footnote{We can also do a midpoint approximation to find $F_{k,k-1}$ by evaluating the Jacobian at $\left(\hat{x}_{k-1|k-1} + \hat{x}_{k|k-1}\right)/2$.  This should be a much closer approximation to the nonlinear function.  We use this approximation for the Extended Kalman Filter experiments later.}  Numerical accuracy of these methods tends to depend heavily on the nonlinear functions.  If we have linear constraints but a nonlinear $f_{k,k-1}\left(\cdot\right)$ and $h_k\left(\cdot\right)$, we can adapt the Extended Kalman Filter to fit into the framework of the methods described thus far.

\subsection{Nonlinear Equality and Inequality Constraints}

Since equality and inequality constraints we model are often times nonlinear, it is important to make the extension to nonlinear equality and inequality constrained Kalman Filtering for the methods discussed thus far.  Without loss of generality, our discussion here will pertain only to nonlinear inequality constraints.  We can follow the same steps for equality constraints.\footnote{We replace the `$\leq$' sign with an `$=$' sign and the `$\lessapprox$' with an `$\approx$' sign.}  We replace the linear inequality constraint on the state space by the following nonlinear inequality constraint $c\left(x_k\right) = d$, where $c\left(\cdot\right)$ is a vector-valued function.  We can then linearize our constraint, $c\left(x_k\right) = d$, about the current state prediction $\hat{x}_{k|k-1}$, which gives us the following.\footnote{This method is how the Extended Kalman Filter linearizes nonlinear functions for $f_{k,k-1}\left(\cdot\right)$ and $h_k\left(\cdot\right)$.  Here $\hat{x}_{k|k-1}$ can be the state prediction of any of the constrained filters presented thus far and does not necessarily relate to the unconstrained state prediction.}

\begin{equation}
c\left(\hat{x}_{k|k-1}\right) + C \left(x_k - \hat{x}_{k|k-1} \right) \lessapprox d
\end{equation}

Here $C$ is defined as the Jacobian of $c$ evaluated at $\hat{x}_{k|k-1}$.  This indicates then, that the nonlinear constraint we would like to model can be approximated by the following linear constraint

\begin{equation} \label{puenl}
C x_k \lessapprox d + C \hat{x}_{k|k-1} - c\left(\hat{x}_{k|k-1}\right)
\end{equation}

This constraint can be written as $\tilde{C} x_k \leq \tilde{d}$, which is an approximation to the nonlinear inequality constraint.  It is now in a form that can be used by the methods described thus far.

The nonlinearities in both the constraints and the models, $f_{k,k-1}\left(\cdot\right)$ and $h_k\left(\cdot\right)$, could have been linearized using a number of different methods (e.g., a derivative-free method, a higher order Taylor approximation).  Also an iterative method could be used as in the Iterated Extended Kalman Filter \cite{BLK2001}.

\section{Constraining the State Prediction} \label{sec::csp}

We haven't yet discussed whether the state prediction (Equation \eqref{kfsp}) also should be constrained.   Forcing the constraints should provide a better prediction (which is used for forecasting in the Kalman Filter).  Ideally, the transition matrix $F_{k,k-1}$ will take an updated state estimate satisfying the constraints at time $k-1$ and make a prediction that will satisfy the constraints at time $k$.  Of course this may not be the case.  In fact, the constraints may depend on the updated state estimate, which would be the case for nonlinear constraints.  On the downside, constraining the state prediction increases computational cost per iteration.

We propose three methods for dealing with the problem of constraining the state prediction.  The first method is to project the matrix $F_{k,k-1}$ onto the constrained space.  This is only possible for the equality constraints, as there is no trivial way to project $F_{k,k-1}$ to an inequality constrained space.  We can use the same projector as in Equation \eqref{Peq} so we have the following.\footnote{In these three methods, the symmetric weighting matrix $W_k$ can be different.  The resulting $\Upsilon$ can consequently also be different.}  

\begin{equation}
F_{k,k-1}^P = \left(\I - \Upsilon A \right) F_{k,k-1}
\end{equation}

Under the assumption that we have constrained our updated state estimate, this new transition matrix will make a prediction that will keep the estimate in the equality constrained space.  Alternatively, if we weaken this assumption, i.e., we are not constraining the updated state estimate, we could solve the minimization problem below (analogous to Equation \eqref{eq-proj-problem}).  We can also incorporate inequality constraints now.

\begin{equation}
\begin{aligned}
\check{x}_{k|k-1}^{P} =  \argmin_{x} &\  \left(x - \hat{x}_{k|k-1} \right)' W_k \left(x - \hat{x}_{k|k-1} \right) \\
\textnormal{s.t. } & A x = b \\
& C x \leq d
\end{aligned}
\end{equation}





We can constrain the covariance matrix here also, in a similar fashion to the method described in Section~\ref{sec::pue-ineq}.  The third method is to add to the constrained problem the additional constraints below, which ensure that the chosen estimate  will produce a prediction at the next iteration that is also constrained.

\begin{equation}
\begin{aligned}
A_{k+1} F_{k+1,k} x_k &= b_{k+1} \\
C_{k+1} F_{k+1,k} x_k &\leq d_{k+1}
\end{aligned}
\end{equation}

If $A_{k+1}, b_{k+1}, C_{k+1}$ or $d_{k+1}$ depend on the estimate (e.g., if we are linearizing nonlinear functions $a\left(\cdot\right)$ or $b\left(\cdot\right)$, we can use an iterative method, which would resolve $A_{k+1}$ and $b_{k+1}$ using the current best updated state estimate (or prediction), re-calculate the best estimate using $A_{k+1}$ and $b_{k+1}$, and so forth until we are satisfied with the convergence.  This method would be preferred since it looks ahead one time-step to choose a better estimate for the current iteration.\footnote{Further, we can add constraints for some arbitrary $n$ time-steps ahead.}  However, it can be far more expensive computationally.

\section{Experiments}

We provide two related experiments here.  We have a car driving along a straight road with thickness 2 meters.  The driver of the car traces a noisy sine curve (with the noise lying only in the frequency domain).   The car is tagged with a device that transmits the location within some known error.  We would like to track the position of the car.  In the first experiment, we filter over the noisy data with the knowledge that the underlying function is a noisy sine curve.  The inequality constrained methods will constrain the estimates to only take values in the interval $[-1,1]$.  In the second experiment, we do not use the knowledge that the underlying curve is a sine curve.  Instead we attempt to recover the true data using an autoregressive model of order 6 \cite{BJ1976}.  We do, however, assume our unknown function only takes values in the interval $[-1,1]$, and we can again enforce these constraints when using the inequality constrained filter.

The driver's path is generated using the nonlinear stochastic process given by Equation \eqref{kfsm-nl}.  We start with the following initial point.

\begin{equation} \label{ickf1-x0}
x_0 = \begin{bmatrix}
		0 \text{\ m}\\
		0 \text{\ m}
	\end{bmatrix}
\end{equation}

Our vector-valued transition function will depend on a discretization parameter $T$ and can be expressed as below.  Here, we choose $T$ to be $\pi/10$, and we run the experiment from an initial time of 0 to a final time of $10 \pi$.

\begin{equation}
f_{k,k-1} = \begin{bmatrix}
		\left(x_{k-1}\right)_1 + T \\
		\sin \left(\left(x_{k-1}\right)_1 + T \right) 
	\end{bmatrix}
\end{equation}

And for the process noise we choose the following.

\begin{equation}
Q_{k,k-1} = \begin{bmatrix}
	0.1 \text{\ m}^2 & 0 \\
	0 & 0 \text{\ m}^2
	\end{bmatrix}
\end{equation}

The driver's path is drawn out by the second element of the vector $x_k$ -- the first element acts as an underlying state to generate the second element, which also allows a natural method to add noise in the frequency domain of the sine curve while keeping the process recursively generated.

\subsection{First Experiment}

\begin{figure}[h!]  
\begin{center}
\includegraphics[angle=270,width=\columnwidth]{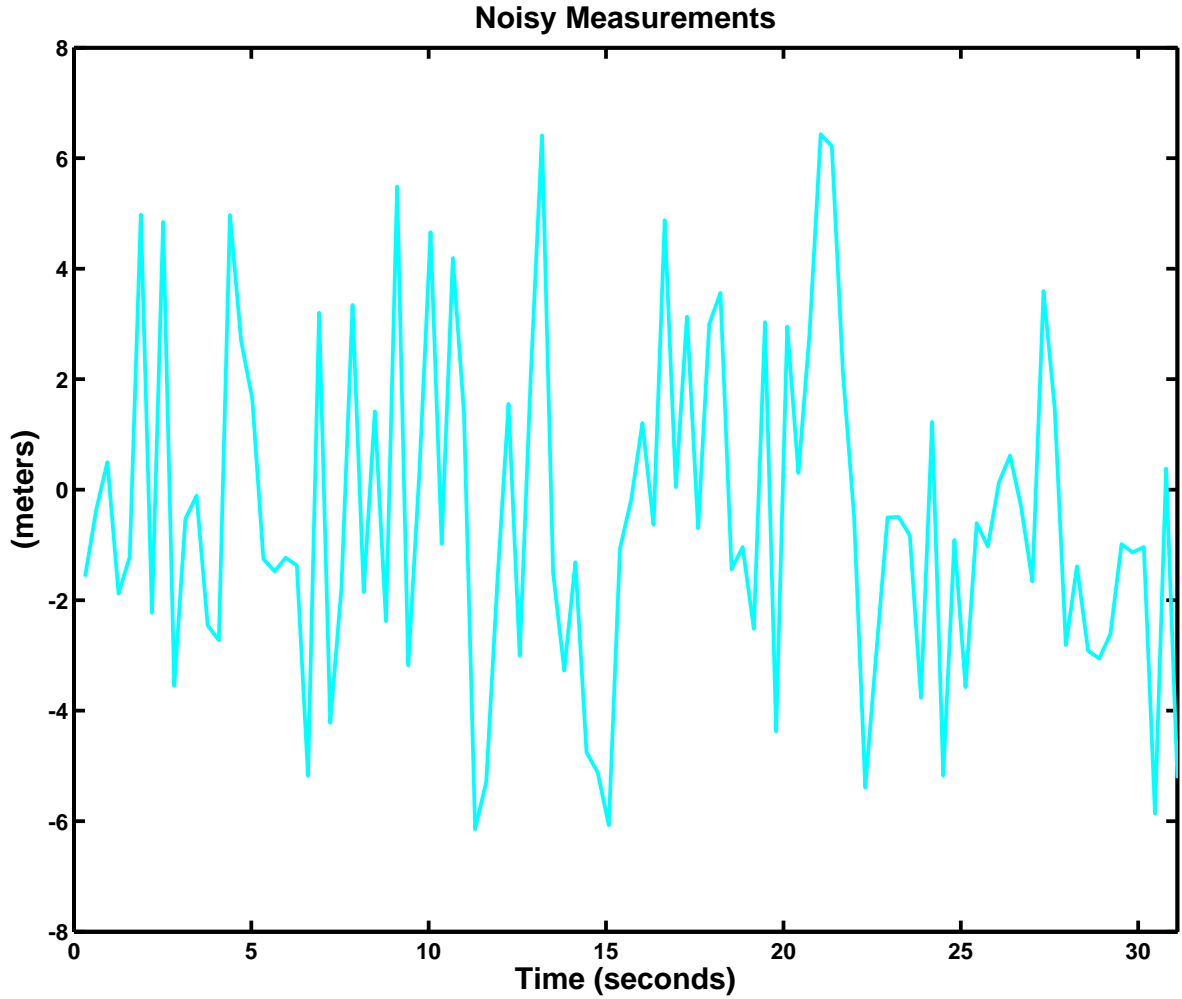}
\end{center}
\caption[Noisy Measurements - Sine Curve]{We take our sine curve, which is already noisy in the frequency domain (due to process noise), and add measurement noise.  The underlying sine curve is significantly masked.}
\label{fig-ickf1}
\end{figure}

To create the measurements, we use the model from Equation \eqref{kfmm}, where $H_k$ is the square identity matrix of dimension 2.  We choose $R_k$ as below to noise the data.  This considerably masks the true underlying data as can be seen in Fig.~\ref{fig-ickf1}.\footnote{The figure only shows the noisy sine curve, which is the second element of the measurement vector.  The first element, which is a noisy straight line, isn't plotted.}

\begin{equation} \label{ickf1-R}
R_{k} = \begin{bmatrix}
	10 \text{\ m}^2 & 0 \\
	0 & 10 \text{\ m}^2
	\end{bmatrix}
\end{equation}

For the initial point of our filters, we choose the following point, which is different from the true initial point given in Equation \eqref{ickf1-x0}.

\begin{equation}
\hat{x}_{0|0} = \begin{bmatrix}
		0 \text{\ m}\\
		1 \text{\ m}
	\end{bmatrix}
\end{equation}

Our initial covariance is given as below.\footnote{Nonzero off-diagonal elements in the initial covariance matrix often help the filter converge more quickly}.  

\begin{equation}
P_{0|0} = \begin{bmatrix}
	1 \text{\ m}^2 & 0.1\\
	0.1 & 1 \text{\ m}^2
	\end{bmatrix}
\end{equation}

In the filtering, we use the information that the underlying function is a sine curve, and our transition function $f_{k,k-1}$ changes to reflect a recursion in the second element of $x_k$ -- now we will add on discretized pieces of a sine curve to our previous estimate.  The function is given explicitly below.

\begin{equation}
f_{k,k-1} = \begin{bmatrix}
		\left(x_{k-1}\right)_1 + T \\
		\left(x_{k-1}\right)_1 + \sin \left(\left(x_{k-1}\right)_1 + T \right) - \sin \left(\left(x_{k-1}\right)_1\right) 
	\end{bmatrix}
\end{equation}

For the Extended Kalman Filter formulation, we will also require the Jacobian of this matrix denoted $F_{k,k-1}$, which is given below.

\begin{equation}
F_{k,k-1} = \begin{bmatrix}
	1 & 0 \\
	\cos \left(\left(x_{k-1}\right)_1 + T \right) - \cos \left(\left(x_{k-1}\right)_1\right)  & 1
	\end{bmatrix}
\end{equation}

The process noise $Q_{k,k-1}$, given below, is chosen similar to the noise used in generating the simulation, but is slightly larger to encompass both the noise in our above model and to prevent divergence due to numerical roundoff errors.  The measurement noise $R_k$ is chosen the same as in Equation \eqref{ickf1-R}.

\begin{equation}
Q_{k,k-1} = \begin{bmatrix}
	0.1 \text{\ m}^2 & 0 \\
	0 & 0.1 \text{\ m}^2
	\end{bmatrix}
\end{equation}

The inequality constraints we enforce can be expressed using the notation throughout the chapter, with $C$ and $d$ as given below.

\begin{equation}
C = \begin{bmatrix}
	0 & 1  \\
	0 & -1 
	\end{bmatrix}
\end{equation}

\begin{equation}
d = \begin{bmatrix}
		1\\
		1
	\end{bmatrix}
\end{equation}

\begin{figure}[h!]  
\begin{center}
\includegraphics[angle=270,width=\columnwidth]{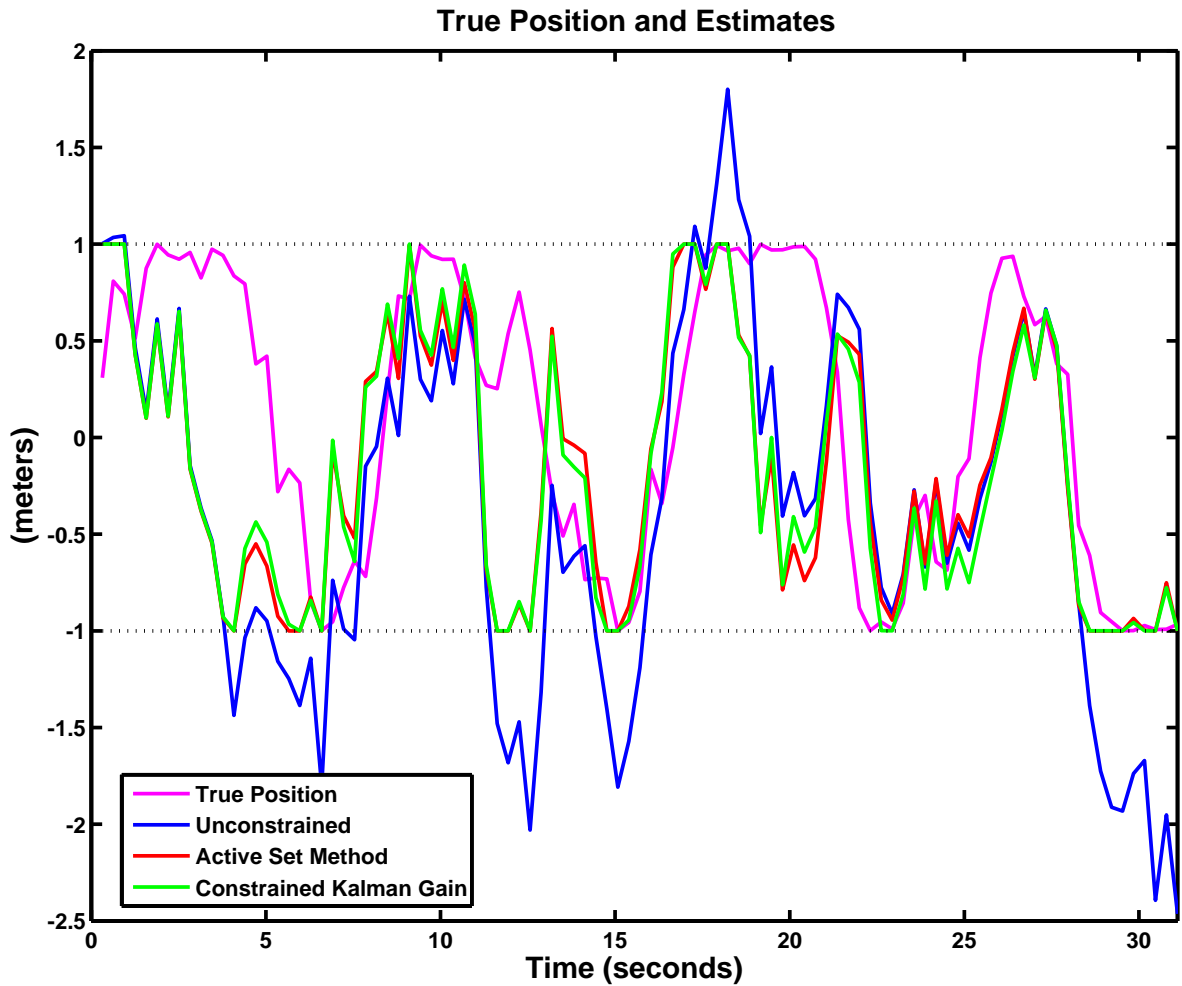}
\end{center}
\caption[True Position and Estimates - Sine Curve]{We show our true underlying state, which is a sine curve noised in the frequency domain, along with the estimates from the unconstrained Kalman Filter, and both of our inequality constrained modifications.  We also plotted dotted horizontal lines at the values -1 and 1.  Both inequality constrained methods do not allow the estimate to leave the constrained space.}
\label{fig-ickf2}
\end{figure}

These constraints force the second element of the estimate $x_{k|k}$ (the sine portion) to lie in the interval $[-1,1]$.  We do not have any equality constraints in this experiment.  We run the unconstrained Kalman Filter and both of the constrained methods discussed previously.  A plot of the true position and estimates is given in Fig.~\ref{fig-ickf2}.  Notice that both constrained methods force the estimate to lie within the constrained space, while the unconstrained method can violate the constraints.

\subsection{Second Experiment}

In the previous experiment, we used the knowledge that the underlying function was a noisy sine curve.  If this is not known, we face a significantly harder estimation problem.  Let us assume nothing about the underlying function except that it must take values in the interval $[-1,1]$.  A good model for estimating such an unknown function could be an autoregressive model.  We can compare the unconstrained filter to the two constrained methods again using these assumption and an autoregressive model of order 6, or AR(6) as it is more commonly referred to. 

In the previous example, we used a large measurement noise $R_k$ to emphasize the gain achieved by using the constraint information.  Such a large $R_k$ is probably not very realistic, and when using an autoregressive model, it will be hard to track such a noisy signal.  To generate the measurements, we again use Equation \eqref{kfmm}, this time with $H_k$ and $R_k$ as given below.

\begin{equation}
H_k = \begin{bmatrix}
		0 & 1
	\end{bmatrix}
\end{equation}

\begin{equation}
R_k = \begin{bmatrix}
		0.5 \text{\ m}^2
	\end{bmatrix}
\end{equation}

Our state will now be defined using the following 13-vector, in which the first element is the current estimate, the next five elements are lags, the six elements afterwards are coefficients on the current estimate and the lags, and the last element is a constant term.

\begin{equation}
\hat{x}_{k|k} = \begin{bmatrix}
		y_k  & y_{k-1} & \cdots &  y_{k-5} & \alpha_1 & \alpha_2 & \cdots & \alpha_7
	\end{bmatrix}'
\end{equation}

Our matrix $H_k$ in the filter is a row vector with the first element 1, and all the rest as 0, so $y_{k|k-1}$ is actually our prediction $\hat{z}_{k|k-1}$ in the filter, describing where we believe the expected value of the next point in the time-series to lie.  For the initial state, we choose a vector of all zeros, except the first and seventh element, which we choose as 1.  This choice for the initial conditions leads to the first prediction on the time series being 1, which is incorrect as the true underlying state has expectation 0.  For the initial covariance, we choose $\I_{\by{13}{13}}$ and add $0.1$ to all the off-diagonal elements.\footnote{The bracket subscript notation is used through the remainder of this paper to indicate the size of zero matrices and identity matrices.}  The transition function $f_{k,k-1}$ for the AR(6) model is given below.

\begin{equation}
\begin{bmatrix}
		\min\left(1, \max\left(-1, \alpha_1 y_{k-1} + \cdots + \alpha_6 y_{k-6} + \alpha_7 \right) \right)\\
		\min\left(1, \max\left(-1,y_{k-1} \right) \right)\\
		\min\left(1, \max\left(-1,y_{k-2} \right) \right)\\
		\min\left(1, \max\left(-1,y_{k-3} \right) \right)\\
		\min\left(1, \max\left(-1,y_{k-4} \right) \right)\\
		\min\left(1, \max\left(-1,y_{k-5} \right) \right)\\
		\alpha_1 \\
		\alpha_2 \\
		\vdots \\
		\alpha_6 \\
		\alpha_7
	\end{bmatrix}
\end{equation}

Putting this into recursive notation, we have the following.


\begin{equation}
\begin{bmatrix}
		\min\left(1, \max\left(-1, \left(x_{k-1}\right)_7 \left(x_{k-1}\right)_1 + \cdots + \left(x_{k-1}\right)_{13} \right) \right)\\
		\min\left(1, \max\left(-1, \left(x_{k-1}\right)_1 \right) \right)\\
		\min\left(1, \max\left(-1, \left(x_{k-1}\right)_2 \right) \right)\\
		\min\left(1, \max\left(-1, \left(x_{k-1}\right)_3 \right) \right)\\
		\min\left(1, \max\left(-1, \left(x_{k-1}\right)_4 \right) \right)\\
		\min\left(1, \max\left(-1, \left(x_{k-1}\right)_5 \right) \right)\\
		\left(x_{k-1}\right)_7 \\
		\left(x_{k-1}\right)_8 \\
		\vdots \\
		\left(x_{k-1}\right)_{12} \\
		\left(x_{k-1}\right)_{13}
	\end{bmatrix}
\end{equation}

The Jacobian of $f_{k,k-1}$ is given below.  We ignore the $\min \left( \cdot \right)$ and $\max \left( \cdot \right)$ operators since the derivative is not continuous across them, and we can reach the bounds by numerical error.  Further, when enforced, the derivative would be 0, so by ignoring them, we are allowing our covariance matrix to be larger than necessary as well as more numerically stable.

\begin{equation}
\begin{bmatrix}
\begin{BMAT}{c.c}{c.c}
	\begin{BMAT}{c.c}{c.c}
	      	\begin{BMAT}{cc}{c} 
			\left(x_{k-1}\right)_7 & \cdots   
		\end{BMAT} &  \left(x_{k-1}\right)_{12} \\
		\I_{\by{5}{5}} & 0_{\by{5}{1}}
	\end{BMAT} & \begin{BMAT}{c}{c.c}
		\begin{BMAT}{cccc}{c}
			\left(x_{k-1}\right)_{1} & \cdots & \left(x_{k-1}\right)_{6} & 1  \\
		\end{BMAT} \\
		0_{\by{5}{7}}
	\end{BMAT} \\
	0_{\by{7}{6}} & \I_{\by{7}{7}}
\end{BMAT}
\end{bmatrix}
\end{equation}

For the process noise, we choose $Q_{k,k-1}$ to be a diagonal matrix with the first entry as 0.1 and all remaining entries as $10^{-6}$ since we know the prediction phase of the autoregressive model very well.  The inequality constraints we enforce can be expressed using the notation throughout the chapter, with $C$ as given below and $d$ as a 12-vector of ones.

\begin{equation}
C = \begin{bmatrix}
\begin{BMAT}{c.c}{c}
	\begin{BMAT}{c}{c.c}
		\I_{\by{6}{6}} \\
		-\I_{\by{6}{6}}
	\end{BMAT} & 0_{\by{12}{7}}
\end{BMAT}
\end{bmatrix}
\end{equation}

\begin{figure}[h!]  
\begin{center}
\includegraphics[angle=270,width=\columnwidth]{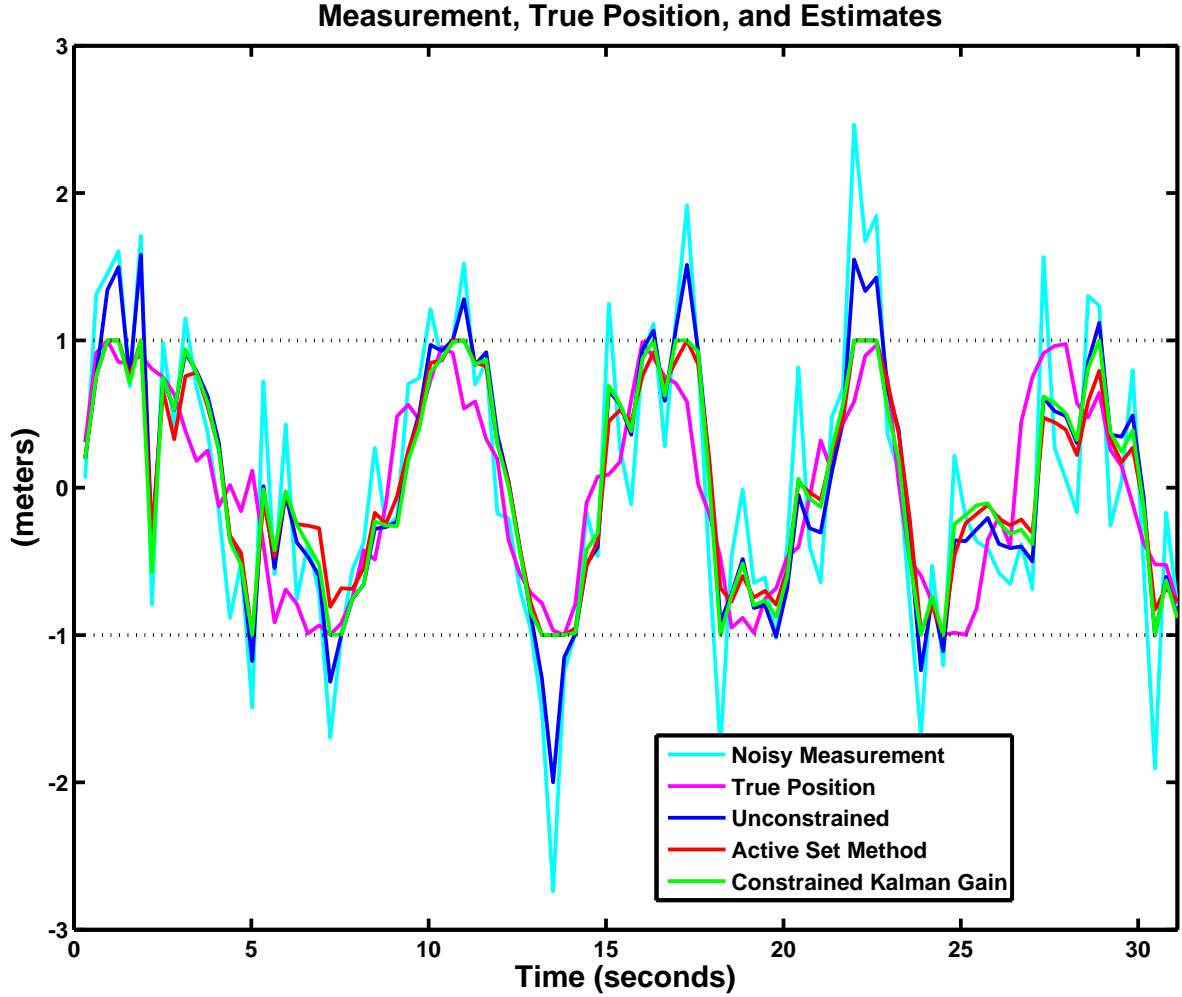}
\end{center}
\caption[Measurement, True Position, and Estimates - Sine Curve and AR(6)]{We show our true underlying state, which is a sine curve noised in the frequency domain, the noised measurements, and the estimates from the unconstrained and both inequality constrained filters.  We also plotted dotted horizontal lines at the values -1 and 1.  Both inequality constrained methods do not allow the estimate to leave the constrained space.}
\label{fig-ickfb}
\end{figure}

These constraints force the current estimate and all of the lags to take values in the range $[-1,1]$.  As an added feature of this filter, we are also estimating the lags at each iteration using more information although we don't use it -- this is a fixed interval smoothing.  In Fig.~\ref{fig-ickfb}, we plot the noisy measurements, true underlying state, and the filter estimates.  Notice again that the constrained methods keep the estimates in the constrained space.  Visually, we can see the improvement particularly near the edges of the constrained space.

\section{Conclusions}

We've provided two different formulations for including constraints into a Kalman Filter.  In the equality constrained framework, these formulations have analytic formulas, one of which is a special case of the other.  In the inequality constrained case, we've shown two numerical methods for constraining the estimate.  We also discussed how to constrain the state prediction and how to handle nonlinearities.  Our two examples show that these methods ensure the estimate lies in the constrained space, which provides a better estimate structure.

\appendices

\section{Kron and Vec} \label{app::kv}

In this appendix, we provide some definitions used earlier in the chapter.  Given matrix $A \in \mathbb{R}^{ m \times n}$ and $B \in \mathbb{R}^{p \times q}$, we can define the right Kronecker product as below.\footnote{The indices $m,n,p$, and $q$ and all matrix definitions are independent of any used earlier.  Also, the subscript notation $a_{1,n}$ denotes the element in the first row and $n$-th column of $A$, and so forth.}

\begin{equation}
\left( A \otimes B \right) = \begin{bmatrix}
a_{1,1} B & \cdots & a_{1,n} B \\
\vdots & \ddots & \vdots \\
a_{m,1} B & \cdots & a_{m,n} B
\end{bmatrix}
\end{equation}

Given appropriately sized matrices $A, B, C,$ and $D$ such that all operations below are well-defined, we have the following equalities.

\begin{equation} \label{kron-trans}
\left( A \otimes B \right)' = \left( A' \otimes B' \right)
\end{equation}

\begin{equation} \label{kron-inv}
\left( A \otimes B \right) ^{-1} = \left( A^{-1} \otimes B^{-1} \right)
\end{equation}

\begin{equation} \label{kron-dist}
\left( A \otimes B \right) \left( C \otimes D \right) = \left( AC \otimes BD \right)
\end{equation}

We can also define the vectorization of an $\by{m}{n}$ matrix $A$, which is a linear transformation on a matrix that stacks the columns iteratively to form a long vector of size $\by{mn}{1}$, as below.

\begin{equation}
\veco{A} = \begin{bmatrix}
a_{1,1} \\
\vdots \\
a_{m,1} \\
a_{1,2} \\
\vdots \\
a_{m,2} \\
\vdots \\
a_{1,n} \\
\vdots \\
a_{m,n}
\end{bmatrix}
\end{equation}

Using the vec operator, we can state the trivial definition below.

\begin{equation} \label{vec-sum}
\veco{A+B} = \veco{A} + \veco{B}
\end{equation}

Combining the vec operator with the Kronecker product, we have the following.

\begin{equation} \label{vec-ab}
\veco{AB} = \kron{B'}{\I} \veco{A}
\end{equation}

\begin{equation} \label{vec-abc}
\veco{ABC} = \left(C' \otimes A \right) \veco{B}
\end{equation}

We can express the trace of a product of matrices as below.

\begin{equation} \label{tr-ab}
\tr{AB} = \veco{B'}'\veco{A}
\end{equation}

\begin{subequations}
\begin{align}
\tr{ABC} &= 
	\label{trace-1} \veco{B}' \left(\I \otimes C\right) \veco{A} \\
&= 
	\label{trace-2} \veco{A}' \left(\I \otimes B \right) \veco{C} \\
&=
	\label{trace-3} \veco{A}' \left(C \otimes \I \right) \veco{B}
\end{align}
\end{subequations}

For more information, please see \cite{LT1985}.
\section{Analytic Block Representation for the inverse of a Saddle Point Matrix} \label{app::spm}

$M_S$ is a saddle point matrix if it has the block form below.\footnote{The subscript $S$ notation is used to differentiate these matrices from any matrices defined earlier.}

\begin{equation} \label{spm}
M_S =
	\begin{bmatrix}
		A_S & B_S' \\
		B_S & -C_S
	\end{bmatrix}
\end{equation}

In the case that $A_S$ is nonsingular and the Schur complement $J_S = -\left(C_S + B_S A_S^{-1} B_S'\right)$ is also nonsingular in the above equation, it is known that the inverse of this saddle point matrix can be expressed analytically by the following equation (see e.g., \cite{BGL2005}).

\begin{equation}
M_S^{-1} =
	\begin{bmatrix}
		A_S^{-1} + A_S^{-1} B_S'  J_S^{-1} B_S A_S^{-1} & -A_S^{-1} B_S' J_S^{-1} \\
		-J_S^{-1} B_S A_S^{-1} & J_S^{-1}
	\end{bmatrix}
\end{equation}

\section{Solution to the system $Mn=p$} \label{app::Mnp}

Here we solve the system $Mn=p$ from Equations \eqref{M-matrix}, \eqref{n-vector}, and \eqref{p-vector}, re-stated below, for vector $n$.

\begin{equation} \label{Mnp}
\begin{bmatrix}
	2  \kronnp{S_k}{\I} &  \nu_{k} \otimes A' \\
	 \nu_{k}' \otimes A  & 0_{\by{q}{q}}
\end{bmatrix} \begin{bmatrix}
	l \\
	\lambda
\end{bmatrix} = \begin{bmatrix}
	0_{\by{mn}{1}} \\
	b - A \hat{x}_{k|k}
\end{bmatrix}
\end{equation}

$M$ is a saddle point matrix with the following equations to fit the block structure of Equation \eqref{spm}.\footnote{We use Equation \eqref{kron-trans} with $B_S'$ to arrive at the same term for $B_s$ in Equation \eqref{Mnp}.}

\begin{align}
A_S & =  2  \kronnp{S_k}{\I} \\
B_S & = \nu_{k}' \otimes A  \\
C_S & = 0_{\by{q}{q}}
\end{align}

We can calculate the term $A_S^{-1} B_S'$.

\begin{subequations}
\begin{align}
A_S^{-1} B_S' & = \left[ 2\kron{S_k}{\I}\right]^{-1} \left( \nu_{k}' \otimes A \right)'  \\
&\stackrel{\eqref{kron-trans}\eqref{kron-inv}}{=}  \frac{1}{2} \kron{S_k^{-1}}{\I} \left( \nu_{k} \otimes A' \right) \\
&\stackrel{\eqref{kron-dist}}{=}  \frac{1}{2} \left( S_k^{-1} \nu_k \right) \otimes A'
\end{align}
\end{subequations}

And as a result we have the following for $J_S$.

\begin{subequations}
\begin{align}
J_S & = - \frac{1}{2} \left( \nu_{k}' \otimes A \right) \left[ \left( S_k^{-1} \nu_k \right) \otimes A' \right] \\
&\stackrel{\eqref{kron-dist}}{=} - \frac{1}{2} \left( \nu_{k}' S_k^{-1} \nu_k \right) \otimes \left(A A' \right) 
\end{align}
\end{subequations}

$J_S^{-1}$ is then, as below.

\begin{subequations}
\begin{align}
J_S^{-1} & = -2 \left[ \left( \nu_{k}' S_k^{-1} \nu_k \right) \otimes \left( A A' \right)\right]^{-1} \\
&\stackrel{\eqref{kron-inv}}{=} -2 \left(\nu_{k}' S_k^{-1} \nu_k  \right)^{-1} \otimes \left(A A' \right)^{-1}
\end{align}
\end{subequations}

For the upper right block of $M^{-1}$, we then have the following expression.

\begin{subequations}
\begin{align}
A_S^{-1} B_S' J_S^{-1} &= \left[\left( S_k^{-1} \nu_k \right) \otimes A' \right] \left[\left(\nu_{k}' S_k^{-1} \nu_k \right)^{-1} \otimes \left(A A' \right)^{-1}\right] \\
&\stackrel{\eqref{kron-dist}}{=}  \left[S_k^{-1} \nu_k \left(\nu_{k}' S_k^{-1} \nu_k \right)^{-1}\right] \otimes \left[A' \left(A A' \right)^{-1} \right]
\end{align}
\end{subequations}

Since the first block element of $p$ is a vector of zeros, we can solve for $n$ to arrive at the following solution for $l$.

\begin{equation}
\left(\left[S_k^{-1} \nu_k \left(\nu_{k}' S_k^{-1} \nu_k \right)^{-1}\right] \otimes \left[A' \left(A A' \right)^{-1} \right]\right) \left(b - A \hat{x}_{k|k}\right) \\
\end{equation}

The vector of Lagrange Multipliers $\lambda$ is given below.

\begin{equation}
-2 \left[\left(\nu_{k}' S_k^{-1} \nu_k \right)^{-1} \otimes \left(A A' \right)^{-1} \right] \left(b - A \hat{x}_{k|k}\right) 
\end{equation}


\bibliography{ieee-tac3}

\end{document}